\documentclass[dvipdfmx]{article}

\usepackage{amsmath,amssymb,amsthm,mathrsfs}
\usepackage[dvipdfmx]{graphicx}
\newtheorem{theo}{Theorem}

\title{The special values of $L$-functions at $s=1$ of theta products of weight $3$}
\author{Ryojun Ito\thanks{
Department of Mathematics and Informatics, Graduate School of Science, Chiba University, 
Yayoicho 1-33, Inage, Chiba, 263-8522 Japan. E-mail: afua9032@chiba-u.jp,  
2010 Mathematics Subject Classification: 11F27, 11F67, 33C20.       
keywords: theta series, $L$-value for theta products, generalized hypergeometric function. }}

\date{}

\begin{document}

\maketitle

\begin{abstract}
In this paper, we compute the special values of $L$-functions at $s=1$ of some theta products of weight $3$, 
and express them in terms of special values of generalized hypergeometric functions.
\end{abstract}

\section{Introduction and Main Results}
For a modular form $f$ of weight $k$ with the $q$-expansion $f(q) = \sum_{n=1}^{\infty} a_{n}q^{n}$, 
we define its $L$-function $L(f, s)$ by 
\begin{align*}
L(f, s) := \sum_{n=1}^{\infty} \frac{a_{n}}{n^{s}}, \hspace{10mm} {\rm Re}(s) > k+1.
\end{align*}
It is well-known (cf. \cite{shi}) that 
$L(f,s)$ has meromorphic continuation to $\mathbb{C}$, and satisfies a functional equation
under $s \leftrightarrow k-s$. 
In this paper, we consider $L$-functions of modular forms which are products of the Jacobi theta series 
$\theta_{2}(q)$, $\theta_{3}(q)$, $\theta_{4}(q)$ (these are modular forms of weight $1/2$)
or the Borwein theta series $a(q)$, $b(q)$, $c(q)$ (these are modular forms of weight $1$). For the definitions 
of these theta series, see section 2 and 3 respectively.

In 2010s, for some modular forms, it was found that special values of $L$-functions can be expressed in terms of 
special values of generalized hypergeometric functions
\begin{align*}
{_{p+1}F_{p}}\left[ 
\left.
\begin{matrix}
a_{1},a_{2}, \dots , a_{p+1} \\
b_{1}, b_{2}, \dots , b_{p} 
\end{matrix}
\right| z
\right]
:=
\sum_{n=0}^{\infty} \frac{(a_{1})_{n} \cdots (a_{p+1})_{n}}{(b_{1})_{n} \cdots (b_{p})_{n}} \frac{z^{n}}{(1)_{n}},
\end{align*}
where $(a)_{n} := \Gamma(a+n) / \Gamma(a)$ denotes the Pochhammer symbol. For example, 
\begin{enumerate}
\item Otsubo \cite{otsubo} expressed $L(f, 2)$ for some theta products $f(q)$ of weight $2$ 
in terms of ${_{3}F_{2}}(1)$ by using his regulator formula and Bloch's theorem.

\item Rogers \cite{rog}, Rogers-Zudilin \cite{rz}, Zudilin \cite{zud} and 
the author \cite{ito} expressed $L(f, 2)$ for some theta products $f(q)$ of weight $2$ in terms of ${_{3}F_{2}}(1)$ by an 
analytic method. Further, Zudilin \cite{zud} expressed $L(f, 3)$ for the theta product corresponding 
to the elliptic curve of conductor $32$ in terms of ${_{4}F_{3}}(1)$.

\item Rogers-Wan-Zucker \cite{rwz} expressed $L(f, 2)$ (resp. $L(f, 3)$, $L(f,4)$) for some quotients $f(q)$ of 
the Dedekind eta function $\eta(q) := q^{1/24}\prod_{n=1}^{\infty} (1-q^{n})$ of weight $3$ (resp. $4$, $5$) 
in terms of special values of generalized hypergeometric functions or the gamma function $\Gamma(s)$ by an analytic method.

\end{enumerate}
Note that all of these results are in the case of cusp forms.

In this paper, we compute $L(f, 1)$ for theta products (not necessarily cusp forms) of weight $3$, 
and express them in terms of special values of generalized hypergeometric functions. 

Our main results are the following.

\begin{theo}
We have the following table.
\end{theo}
\begin{center}
\begin{tabular}{|c|c|c|} \hline
& $f(q)$ & $L(f,1)$  \\ \hline \hline
& & \\
(i) & $\displaystyle\frac{1}{2}\theta_{2}(q^{4})\theta_{3}^{4}(q^{4})\theta_{4}(q^{8})$ &
$\displaystyle 
\frac{\sqrt{\sqrt{2}+1}\Gamma^{2}\left( \frac{1}{4} \right)}{8\sqrt{\pi}}
{_{3}F_{2}}\left[ 
\left.
\begin{matrix}
\frac{1}{4},\frac{1}{2}, \frac{1}{2} \\
\frac{3}{8}, 1
\end{matrix}
\right| 1
\right]$ \\
& & \\
(ii) &$\displaystyle\frac{1}{2}\theta_{2}(q^{4})\theta_{3}^{3}(q^{4})\theta_{4}^{2}(q^{4})$ &
$\displaystyle 
\frac{\Gamma^{2}\left( \frac{1}{4} \right)}{8\sqrt{2\pi}}
{_{3}F_{2}}\left[ 
\left.
\begin{matrix}
\frac{1}{4},\frac{1}{2}, \frac{1}{2} \\
\frac{3}{4}, 1
\end{matrix}
\right| 1
\right]$ \\
& & \\
(iii) &$\displaystyle\frac{1}{2}\theta_{2}(q^{4})\theta_{3}^{2}(q^{4})\theta_{4}^{3}(q^{4})$ &
$\displaystyle 
\frac{\pi}{4\sqrt{2}}
{_{3}F_{2}}\left[ 
\left.
\begin{matrix}
\frac{1}{4},\frac{1}{2}, \frac{1}{2} \\
1, 1
\end{matrix}
\right| 1
\right]$ \\
& & \\
(iv) &$\displaystyle\frac{1}{2}\theta_{2}(q^{4})\theta_{4}^{5}(q^{4})$ &
$\displaystyle\frac{\Gamma^{2}\left( \frac{1}{4} \right)}{16\sqrt{\pi}}
{_{3}F_{2}}\left[ 
\left.
\begin{matrix}
\frac{1}{4},\frac{1}{2},\frac{1}{2} \\
\frac{3}{2},1
\end{matrix}
\right| 1
\right]$ \\
& & \\
(v) &$\displaystyle\frac{1}{2}\theta_{2}(q^{4})\theta_{4}^{5}(q^{8})$ &
$\displaystyle\frac{\sqrt{\sqrt{2}-1} \Gamma^{4}\left( \frac{1}{4} \right)}{16\pi^{2}} $ \\
& & \\
(vi) &$\displaystyle\frac{1}{4}\theta_{2}^{2}(q^{2})\theta_{3}^{3}(q^{2})\theta_{4}(q^{2})$ & 
$\displaystyle\frac{\Gamma^{2}\left( \frac{1}{4} \right)}{8\sqrt{2\pi}}
{_{3}F_{2}}\left[ 
\left.
\begin{matrix}
\frac{1}{2}, \frac{1}{2}, \frac{1}{2} \\
\frac{3}{4}, 1
\end{matrix}
\right| 1
\right]$ \\
& & \\
(vii) &$\displaystyle\frac{1}{4}\theta_{2}^{2}(q^{2})\theta_{3}(q^{2})\theta_{4}^{3}(q^{2})$ &
$\displaystyle 
\frac{\Gamma^{2}\left( \frac{3}{4} \right)}{2\sqrt{2\pi}}
{_{3}F_{2}}\left[ 
\left.
\begin{matrix}
\frac{1}{2},\frac{1}{2}, \frac{1}{2} \\
\frac{5}{4}, 1
\end{matrix}
\right| 1
\right]$ \\
& & \\
(viii) &$\displaystyle\frac{1}{4}\theta_{2}^{2}(q^{2})\theta_{4}^{4}(q^{2})$ &
$\displaystyle 
\frac{1}{4}
{_{3}F_{2}}\left[ 
\left.
\begin{matrix}
\frac{1}{2},\frac{1}{2}, \frac{1}{2} \\
\frac{3}{2}, 1
\end{matrix}
\right| 1
\right]$ \\
& & \\
(ix) &$\displaystyle\frac{1}{8}\theta_{2}^{3}(q^{4})\theta_{3}^{2}(q^{4})\theta_{4}(q^{4})$ &
$\displaystyle 
\frac{\pi}{16\sqrt{2}}
{_{3}F_{2}}\left[ 
\left.
\begin{matrix}
\frac{3}{4},\frac{1}{2}, \frac{1}{2} \\
1, 1
\end{matrix}
\right| 1
\right]$ \\
& & \\
(x) &$\displaystyle\frac{1}{8}\theta_{2}^{3}(q^{4})\theta_{3}(q^{4})\theta_{4}^{2}(q^{4})$ &
$\displaystyle 
\frac{\Gamma^{2}\left( \frac{3}{4} \right)}{8\sqrt{2\pi}}
{_{3}F_{2}}\left[ 
\left.
\begin{matrix}
\frac{3}{4},\frac{1}{2}, \frac{1}{2} \\
\frac{5}{4}, 1
\end{matrix}
\right| 1
\right]$ \\
& & \\
(xi) &$\displaystyle\frac{1}{8}\theta_{2}^{3}(q^{4})\theta_{4}^{3}(q^{4})$ &
$\displaystyle 
\frac{\pi}{32\sqrt{2}}$ \\
& & \\
\hline
\end{tabular}
\end{center}
\begin{center}
\begin{tabular}{|c|c|c|}
\hline 
& & \\
(xii) &$\displaystyle\frac{1}{8}\theta_{2}^{3}(q^{4})\theta_{4}^{3}(q^{8})$ &
$\displaystyle 
\frac{\Gamma\left( \frac{1}{8} \right)\Gamma\left( \frac{1}{4} \right)\Gamma^{3}\left( \frac{3}{8} \right)}
{128\pi^{5/2}}$ \\
& & \\
(xiii) &$\displaystyle\frac{1}{16}\theta_{2}^{4}(q)\theta_{4}^{2}(q)$ & 
$\displaystyle\frac{\pi}{16}$ \\ 
& & \\
(xiv) &$\displaystyle\frac{1}{32}\theta_{2}^{5}(q^{4})\theta_{4}(q^{4})$ & 
$\displaystyle\frac{\Gamma^{2}\left( \frac{1}{4} \right)}{256\sqrt{\pi}}
{_{3}F_{2}}\left[ 
\left.
\begin{matrix}
\frac{5}{4},\frac{1}{2},\frac{1}{2} \\
\frac{3}{2},1
\end{matrix}
\right| 1
\right]$ \\ 
& & \\
(xv) &$\displaystyle\frac{1}{32}\theta_{2}^{5}(q^{4})\theta_{4}(q^{8})$ & 
$\displaystyle\frac{\sqrt{\sqrt{2}+1}\Gamma^{4}\left( \frac{1}{4} \right)}{128\pi^{2}}$ \\ 
& & \\ \hline
\end{tabular}

\end{center}

\begin{theo}
We have the following table.
\end{theo}

\begin{center}
\begin{tabular}{|c|c|c|} \hline
& $f(q)$ & $L(f,1)$  \\ \hline \hline
& & \\
(i) &$\displaystyle \frac{1}{3}a(q^{3})c(q^{3})b(q^{3})$ & 
$\displaystyle\frac{\Gamma^{6}\left( \frac{1}{3}\right)}{8\sqrt{3}\pi^{3}}$ \\
& & \\
(ii) &$\displaystyle \frac{1}{3}c(q^{3})b^{2}(q^{3})$ & 
$\displaystyle\frac{2\pi}{9\sqrt{3}}
{_{3}F_{2}}\left[ 
\left.
\begin{matrix}
\frac{1}{3}, \frac{1}{3}, \frac{2}{3} \\
1, 1
\end{matrix}
\right| 1
\right]$ \\
& & \\
(iii) &$\displaystyle \frac{1}{9}c^{2}(q^{3})b(q^{3})$ &
$\displaystyle\frac{2\pi}{27\sqrt{3}}
{_{3}F_{2}}\left[ 
\left.
\begin{matrix}
\frac{1}{3}, \frac{2}{3}, \frac{2}{3} \\
1, 1
\end{matrix}
\right| 1
\right]$ \\
& & \\
(iv) &$\displaystyle\frac{1}{3}c(q^{3})b^{2}(q^{9})$ & 
$\displaystyle\frac{2\pi}{3^{11/6}}
{_{3}F_{2}}\left[ 
\left.
\begin{matrix}
\frac{1}{9}, \frac{4}{9}, \frac{7}{9} \\
1, 1
\end{matrix}
\right|1 \right]$ \\
& & \\
(v) &$\displaystyle\frac{1}{9}c^{2}(q^{3})b(q^{9})$ & 
$\displaystyle\frac{2\pi}{3^{13/6}}
{_{3}F_{2}}\left[ 
\left.
\begin{matrix}
\frac{2}{9}, \frac{5}{9}, \frac{8}{9} \\
1, 1
\end{matrix}
\right|1 \right]$ \\ 
& & \\ \hline
\end{tabular}

\end{center}
Here we normalized $f(q)=\sum_{n=1}^{\infty} a_{n}q^{n}$ so that $a_{1} = 1$.

We remark that, for theta products $f(q)$ of weight $3$ which are considered in \cite[Theorem 5]{rwz}, 
the values $L(f, 2)$ (hence $L(f, 1)$ by the functional equation) 
are expressed in terms of special values of the gamma function, not generalized hypergeometric functions.
It is new that the values of $L$-functions at $s=1$ of theta products of weight $3$ are expressed 
in terms of special values of generalized hypergeometric functions.

Our strategy to compute $L(f,1)$ is the same as that used in \cite{ito}, \cite{rog}, \cite{rwz}, \cite{rz} and \cite{zud}.
For a modular form $f(q)=\sum_{n=1}^{\infty} a_{n}q^{n}$ and $m \in \mathbb{Z}_{\geqq 1}$, 
the value $L(f,m)$ is obtained by the Mellin transformation of $f(q)$
\begin{align*}
L(f, m) = \frac{(-1)^{m-1}}{\Gamma(m)} \int_{0}^{1} f(q) (\log q)^{m-1} \frac{dq}{q}. 
\end{align*}
The case $m=1$ is special since the logarithm in the integral above vanishes: 
\begin{align}
L(f, 1) = \int_{0}^{1} f(q) \frac{dq}{q}. \label{lfunc}
\end{align}
Since the Jacobi theta series and the Borwein theta series have connections 
with hypergeometric functions (see \eqref{jacobitrans} and \eqref{borweintrans} below), 
we can reduce \eqref{lfunc} to an integral of the form 
\begin{align*}
\int_{0}^{1} (\mbox{polynomials in } \alpha^{l} (1-\alpha)^{m}) {_{2}F_{1}}(\alpha^{n}) \frac{d\alpha}{\alpha(1-\alpha)}.
\end{align*}
Then, we obtain a hypergeometric evaluation of $L(f, 1)$ after some computation.

Finally, we remark that one of the results of Rogers-Wan-Zucker \cite{rwz} can be recovered from our results. 
Let $f(q) = \frac{1}{2} \theta_{2}(q^{4})\theta_{3}^{4}(q^{4})\theta_{4}(q^{4})$. By the Jacobi identity (see section 2), 
we have
\begin{align*}
2f(q) = \theta_{2}(q^{4})\theta_{3}^{4}(q^{4})\theta_{4}(q^{4}) = 
\theta_{2}(q^{4})\theta_{4}^{5}(q^{4}) + \theta_{2}^{5}(q^{4})\theta_{4}(q^{4}). 
\end{align*}
Therefore, by Theorem 1 (iv), (xiv), we obtain
\begin{align*}
L(f, 1) = \frac{\Gamma^{2}\left( \frac{1}{4} \right)}{16\sqrt{\pi}} \left( 
{_{3}F_{2}}\left[ 
\left.
\begin{matrix}
\frac{1}{4}, \frac{1}{2}, \frac{1}{2}\\
\frac{3}{2}, 1
\end{matrix}
\right| 1
\right] 
+
{_{3}F_{2}}\left[ 
\left.
\begin{matrix}
\frac{5}{4}, \frac{1}{2}, \frac{1}{2}\\
\frac{3}{2}, 1
\end{matrix}
\right| 1
\right] 
\right).
\end{align*}
If we use the following contiguous relation 
\begin{align*}
(b-a)
{_{3}F_{2}}\left[ 
\left.
\begin{matrix}
a,b,c\\
e,f
\end{matrix}
\right| z
\right]  +
a {_{3}F_{2}}\left[ 
\left.
\begin{matrix}
a+1,b,c\\
e,f
\end{matrix}
\right| z
\right]
= b {_{3}F_{2}}\left[ 
\left.
\begin{matrix}
a,b+1,c\\
e,f
\end{matrix}
\right| z
\right] ,
\end{align*}
then we have
\begin{align*}
L(f, 1) = \frac{\Gamma^{4}\left( \frac{1}{4} \right)}{8\sqrt{2}\pi^{2}}.
\end{align*} 
This is the second last formula in \cite[Theorem 5]{rwz}.

\section{Proof of Theorem 1}

We define the Jacobi theta series $\theta_{2}(q)$, $\theta_{3}(q)$ and $\theta_{4}(q)$ by
\begin{align*}
\theta_{2}(q) &:= \sum_{n \in \mathbb{Z}} q^{(n+1/2)^{2}},  \\
\theta_{3}(q) &:= \sum_{n \in \mathbb{Z}} q^{n^{2}},  \\
\theta_{4}(q) &:= \sum_{n \in \mathbb{Z}} (-1)^{n}q^{n^{2}}.
\end{align*}
These are many relations between $\theta_{2}(q)$, $\theta_{3}(q)$ and $\theta_{4}(q)$. One of the most important relations is 
the Jacobi identity \cite[p.35, (2.1.10)]{borwein1}
\begin{align*}
\theta_{3}^{4}(q) = \theta_{2}^{4}(q) + \theta_{4}^{4}(q).
\end{align*}
Further, the Jacobi theta series have a connection with hypergeometric functions.
Let $\alpha := \theta_{2}^{4} (q)/ \theta_{3}^{4}(q)$. 
Note that we have $1 - \alpha = \theta_{4}^{4}(q) / \theta_{3}^{4}(q)$ by the Jacobi identity. 
Then we have 
\begin{align}
\theta_{3}^{2}(q) =
{_{2}F_{1}}\left[ \left.
\begin{matrix}
\frac{1}{2}, \frac{1}{2} \\
1
\end{matrix}
\right| \alpha \right], \hspace{10mm}
\frac{dq}{q} = \frac{d\alpha}{\alpha(1-\alpha){_{2}F_{1}}^{2}\left[ \left.
\begin{matrix}
\frac{1}{2}, \frac{1}{2} \\
1
\end{matrix}
\right| \alpha \right]}. \label{jacobitrans}
\end{align}
The former is \cite[p.101, Entry 6]{ramanujan3}, and the latter follows from the former and \cite[p.87, Entry 30]{ramanujan2}.

\begin{proof}[Proof of Theorem 1]
Theorem 1 follows from \eqref{jacobitrans}. 
First, we show the formula (xiii).
For simplicity, we denote 
\begin{align*}
\Gamma\left[ 
\begin{matrix}
a_{1}, \dots, a_{p} \\
b_{1}, \dots, b_{q}
\end{matrix}
\right]
:= \frac{\Gamma(a_{1}) \cdots \Gamma(a_{p})}{\Gamma(b_{1}) \cdots \Gamma(b_{q})}.
\end{align*}
By \eqref{lfunc}, we have 
\begin{align*}
L(f, 1) 
= \frac{1}{16} \int_{0}^{1} \theta_{2}^{4}(q)\theta_{4}^{2}(q) \frac{dq}{q} 
= \frac{1}{16} \int_{0}^{1} \left( \frac{\theta_{2}(q)}{\theta_{3}(q)} \right)^{4} 
\left( \frac{\theta_{4}(q)}{\theta_{3}(q)} \right)^{2} \theta_{3}^{6}(q) \frac{dq}{q}.
\end{align*}
If we use \eqref{jacobitrans}, the integral above is equal to 
\begin{align*}
\frac{1}{16}\int_{0}^{1} \alpha (1-\alpha)^{1/2}
{_{2}F_{1}}\left[ \left.
\begin{matrix}
\frac{1}{2}, \frac{1}{2} \\
1
\end{matrix}
\right| \alpha \right]
\frac{d\alpha}{\alpha(1-\alpha)}.
\end{align*}
Since generalized hypergeometric functions have the integral representation \cite[p.108, (4.1.2)]{slater}
\begin{align*}
&{_{p+1}F_{p}}\left[ 
\left.
\begin{matrix}
a_{1},a_{2}, \dots , a_{p+1} \\
b_{1}, b_{2}, \dots , b_{p} 
\end{matrix}
\right| z
\right]  \\
&=
\Gamma \left[ 
\begin{matrix}
b_{1}  \\
a_{1}, b_{1} - a_{1}
\end{matrix}
\right]
\int_{0}^{1} t^{a_{1}}(1-t)^{b_{1}-a_{1}}
{_{p}F_{p-1}}\left[ 
\left.
\begin{matrix}
a_{2}, \dots , a_{p+1} \\
b_{2}, \dots , b_{p} 
\end{matrix}
\right| zt
\right]
\frac{dt}{t(1-t)},
\end{align*}
we have
\begin{align*}
L(f,1)
=\frac{1}{16}\Gamma\left[
\begin{matrix}
1, \frac{1}{2} \\
\frac{3}{2}
\end{matrix}
\right]
{_{3}F_{2}}\left[ 
\left.
\begin{matrix}
1,\frac{1}{2},\frac{1}{2} \\
\frac{3}{2},1
\end{matrix}
\right| 1
\right].
\end{align*}
Note that this ${_{3}F_{2}}$ reduces to a ${_{2}F_{1}}$. 
By Gauss' summation formula \cite[p.28, (1.7.6)]{slater}
\begin{align*}
{_{2}F_{1}}\left[ 
\left.
\begin{matrix}
a,b \\
c 
\end{matrix}
\right| 1
\right] = 
\Gamma\left[
\begin{matrix}
c, c-a-b \\
c-a, c-b
\end{matrix}
\right] ,
\end{align*}
we obtain
\begin{align*}
L(f,1) = \frac{\Gamma^{2}\left( \frac{1}{2} \right)}{16} = \frac{\pi}{16}.
\end{align*}

Next we show the formula (v).
Since we have $\theta_{3}(q)\theta_{4}(q) = \theta_{4}^{2}(q^{2})$ \cite[p.34, (2.1.7 ii)]{borwein1}, 
we obtain
\begin{align*}
L(f, 1) = \frac{1}{8}\int_{0}^{1} \theta_{2}(q) \theta_{3}^{5/2}(q) \theta_{4}^{5/2}(q) \frac{dq}{q}.
\end{align*}
If we use \eqref{jacobitrans} and the integral representation of generalized hypergeometric functions, 
we obtain
\begin{align*}
L(f,1) = \frac{1}{8}
\Gamma\left[
\begin{matrix}
\frac{1}{4}, \frac{5}{8} \\
\frac{7}{8}
\end{matrix}
\right]
{_{3}F_{2}}\left[ 
\left.
\begin{matrix}
\frac{1}{4},\frac{1}{2},\frac{1}{2} \\
\frac{7}{8},1
\end{matrix}
\right| 1
\right].
\end{align*}
By the Watson summation formula \cite[p.54, (2.3.3.13)]{slater}
\begin{align*}
{_{3}F_{2}}\left[ 
\left.
\begin{matrix}
a,b,c\\
\frac{1+a+b}{2}, 2c
\end{matrix}
\right| 1
\right]
= 
\Gamma\left[
\begin{matrix}
\frac{1}{2}, \frac{1+2c}{2}, \frac{1+a+b}{2}, \frac{1-a-b+2c}{2} \\
\frac{1+a}{2}, \frac{1+b}{2}, \frac{1-a+2c}{2}, \frac{1-b+2c}{2}
\end{matrix}
\right] ,
\end{align*}
we have
\begin{align*}
L(f, 1) = \Gamma\left[
\begin{matrix}
\frac{1}{4}, \frac{4}{8}, \frac{1}{2} \\
\frac{7}{8}, \frac{3}{4}, \frac{3}{4}
\end{matrix}
\right].
\end{align*}
Therefore we obtain the formula if we simplify this $\Gamma$-factor 
using the reflection formula and multiplication formula for $\Gamma(s)$
\begin{align*}
&\Gamma(n)\Gamma(1-n) = \frac{\pi}{\sin n\pi}, \\
&\Gamma(nz) = \frac{n^{nz-1/2}}{(2\pi)^{(n-1)/2}} \prod_{k=0}^{n-1} \Gamma\left(z + \frac{k}{n} \right).
\end{align*}

Similar computations lead to the remaining formulas. Note that we use the Watson summation formula 
for the formulas (xi), (xii) and (xv).

\end{proof}

\section{Proof of Theorem 2}
We define the Borwein theta series $a(q)$, $b(q)$ and $c(q)$ \cite{borwein2}, \cite{bbg} by 
\begin{align*}
a(q) &= \sum_{n, m \in\mathbb{Z}} q^{n^{2}+nm +m^{2}},  \\
b(q) &= \sum_{n, m \in \mathbb{Z}} \omega^{n-m}q^{n^{2}+nm +m^{2}}, \\
c(q) &= \sum_{n, m \in \mathbb{Z}} q^{(n+1/3)^{2} + (n+1/3)(m+1/3) + (m+1/3)^{2}},
\end{align*}
where $\omega$ denotes a primitive 3rd root of unity.

There are many relations analogous to those of the Jacobi theta series. For example, J.M. Borwein and P.B. Borwein proved 
\cite{borwein2}
\begin{align*}
a^{3}(q) = b^{3}(q) + c^{3}(q),
\end{align*}
which is a cubic analogue of the Jacobi identity.
Further, like the Jacobi theta series, the Borwein theta series have a connection with hypergeometric functions. Let 
$\alpha := c^{3}(q) / a^{3}(q)$. Note that we have $1-\alpha = b^{3}(q) / a^{3}(q)$ by the cubic identity above. Then we have
\begin{align}
a(q) = 
{_{2}F_{1}}\left[ \left.
\begin{matrix}
\frac{1}{3}, \frac{2}{3} \\
1
\end{matrix}
\right| \alpha \right], \hspace{10mm}
\frac{dq}{q} = \frac{d\alpha}{\alpha(1-\alpha){_{2}F_{1}}^{2}\left[ \left.
\begin{matrix}
\frac{1}{3}, \frac{2}{3} \\
1
\end{matrix}
\right| \alpha \right]}. \label{borweintrans}
\end{align}
The former is \cite[p.97, (2.26)]{ramanujan5}, 
and the latter follows from the former and \cite[p.87, Entry 30]{ramanujan2}).

\begin{proof}[Proof of Theorem 2]
Like the cases of Jacobi products, Theorem 2 follows from \eqref{borweintrans}. 
For example, the formula (i) follows from the following computations.

By \eqref{lfunc}, we have
\begin{align*}
L(f,1) = \frac{1}{9} \int_{0}^{1} a(q)c(q)b(q) \frac{dq}{q}.
\end{align*}
If we use \eqref{borweintrans}, the integral above is equal to 
\begin{align*}
\frac{1}{9}\int_{0}^{1} \alpha^{1/3}(1-\alpha)^{1/3} 
{_{2}F_{1}}\left[ \left.
\begin{matrix}
\frac{1}{3}, \frac{2}{3} \\
1
\end{matrix}
\right| \alpha \right]
\frac{d\alpha}{\alpha(1-\alpha)}.
\end{align*}
By the integral representation of ${_{3}F_{2}}(z)$, we obtain
\begin{align*}
L(f, 1) = \frac{1}{9}
\Gamma\left[
\begin{matrix}
\frac{1}{3}, \frac{1}{3} \\
\frac{2}{3}
\end{matrix}
\right]
{_{3}F_{2}}\left[ \left.
\begin{matrix}
\frac{1}{3},\frac{1}{3}, \frac{2}{3} \\
\frac{2}{3}, 1
\end{matrix}
\right| 1 \right].
\end{align*}
Note that this ${_{3}F_{2}}$ reduces to a ${_{2}F_{1}}$, hence, by Gauss' theorem, we have 
\begin{align*}
L(f, 1) = \frac{1}{9}\Gamma\left[
\begin{matrix}
\frac{1}{3}, \frac{1}{3}, \frac{1}{3} \\
\frac{2}{3}, \frac{2}{3}, \frac{2}{3}
\end{matrix}
\right].
\end{align*}
If we use the reflection formula to simplify the $\Gamma$-factor, we obtain the formula.

Similar computations lead to the formulas (ii) and (iii).

Next we prove the formula (iv).
Since we have $a(q^{3}) = (a(q) + 2b(q))/3 $, $c(q) =(a(q) -b(q))/3$ \cite[Lemma 2.1 (iii) and (2.1)]{bbg}, we obtain
\begin{align*}
b^{3}(q^{3}) = a^{3}(q^{3}) - c^{3}(q^{3}) = \frac{a^{2}(q)b(q)+a(q)b^{2}(q)+b^{3}(q)}{3}.
\end{align*}
Therefore we have
\begin{align*}
L(f,1) &= \frac{1}{9} \int_{0}^{1} c(q)b^{2}(q^{3}) \frac{dq}{q} \\
&=\frac{1}{3^{8/3}} \int_{0}^{1} \alpha^{1/3} ((1-\alpha)^{1/3} +(1-\alpha)^{2/3} + (1-\alpha) )^{2/3}
{_{2}F_{1}}\left[ 
\left.
\begin{matrix}
\frac{1}{3}, \frac{2}{3} \\
1
\end{matrix}
\right| \alpha
\right]
\frac{d\alpha}{\alpha(1-\alpha)}.
\end{align*}
If we substitute $\alpha \mapsto 1- \alpha^{3}$ and $\alpha \mapsto \frac{1-u}{1+2u}$, the integral above is equal to
\begin{align*}
\frac{1}{3^{5/3}} \int_{0}^{1}
\left(\frac{1-u}{1+2u} \right)^{2/3} \left(\frac{3u}{1+2u} \right)^{1/3}
{_{2}F_{1}}\left[ 
\left.
\begin{matrix}
\frac{1}{3}, \frac{2}{3} \\
1
\end{matrix}
\right| 1-\left( \frac{1-u}{1+2u} \right)^{3}
\right] \frac{du}{u(1-u)}.
\end{align*}
We have the cubic transformation \cite[p.96, Theorem 2.3]{ramanujan5}
\begin{align*}
{_{2}F_{1}}\left[ 
\left.
\begin{matrix}
\frac{a}{3},\frac{a+1}{3}\\
\frac{a+1}{2}
\end{matrix}
\right| 1- \left( \frac{1-x}{1+2x} \right)^{3}
\right]
= (1+2x)^{a}
{_{2}F_{1}}\left[ 
\left.
\begin{matrix}
\frac{a}{3},\frac{a+1}{3}\\
\frac{a+5}{6}
\end{matrix}
\right| x^{3}
\right],
\end{align*}
hence we obtain
\begin{align*}
L(f,1) 
= \frac{1}{3^{4/3}}\int_{0}^{1} 
u^{1/3}(1-u)^{2/3} 
{_{2}F_{1}}\left[ 
\left.
\begin{matrix}
\frac{1}{3}, \frac{2}{3} \\
1
\end{matrix}
\right| u^{3}
\right]
\frac{du}{u(1-u)}.
\end{align*}
If we integrate term-by-term, we have
\begin{align*}
L(f, 1)= \frac{1}{3^{4/3}}\sum_{n=0}^{\infty} \frac{\left(\frac{1}{3}\right)_{n}\left( \frac{2}{3} \right)_{n}}{(1)_{n}^{2}}
\Gamma\left[
\begin{matrix}
3n+\frac{1}{3}, \frac{2}{3} \\
3n+1
\end{matrix}
\right] .
\end{align*}
By the multiplication formula, we obtain
\begin{align*}
\Gamma\left(3n+\frac{1}{3} \right) &= \frac{3^{3n-1/6}}{2\pi}\Gamma\left(n + \frac{1}{9} \right) 
\Gamma\left(n + \frac{4}{9} \right) \Gamma\left(n + \frac{7}{9} \right),  \\
\Gamma\left(3n+1 \right) &= \frac{3^{3n+1/2}}{2\pi}\Gamma\left(n + \frac{1}{3} \right) 
\Gamma\left(n + \frac{2}{3} \right) \Gamma\left(n + 1 \right),
\end{align*}
hence we have 
\begin{align*}
L(f,1) &= \frac{1}{9}
\Gamma\left[ 
\begin{matrix}
\frac{1}{9}, \frac{4}{9}, \frac{7}{9}, \frac{2}{3} \\
\frac{1}{3}, \frac{2}{3}
\end{matrix}
\right]
\sum_{n=0}^{\infty} \frac{\left(\frac{1}{3}\right)_{n}\left( \frac{2}{3} \right)_{n}}{(1)_{n}^{2}} \cdot 
\frac{\left(\frac{1}{9}\right)_{n}\left( \frac{4}{9} \right)_{n} \left( \frac{7}{9} \right)_{n}}{\left( 
\frac{1}{3} \right)_{n}\left( \frac{2}{3} \right)_{n}(1)_{n}} \\
&= \frac{1}{9}
\Gamma\left[ 
\begin{matrix}
\frac{1}{9}, \frac{4}{9}, \frac{7}{9} \\
\frac{1}{3}
\end{matrix}
\right]
{_{3}F_{2}}\left[ 
\left.
\begin{matrix}
\frac{1}{9}, \frac{4}{9}, \frac{7}{9} \\
1, 1
\end{matrix}
\right|1 \right].
\end{align*}
We can simplify the $\Gamma$-factor by using the multiplication formula again, then we obtain the formula.

We can prove the formula (v) by similar computations.

\end{proof}

\section*{Acknowledgment}
This paper is based on the author's doctor's thesis at Chiba University.
I would like to thank Noriyuki Otsubo for valuable comments.

\end{document}